\makeatletter
\IfFileExists{./numapde-latex/numapde-packages.sty}{\providecommand*{\input@path}{}\edef\input@path{{./numapde-latex/}\input@path}}{}
\makeatother

\documentclass{numapde-preprint}

\usepackage{numapde-semantic}
\usepackage{numapde-style}
\usepackage{numapde-local}

\IfFileExists{./numapde-bibliography/numapde.bib}{\addbibresource{./numapde-bibliography/numapde.bib}}{\addbibresource{numapde.bib}}

\hypersetup{
	pdftitle={Fenchel Duality and a Separation Theorem on Hadamard Manifolds},
	pdfauthor={Maurício Silva Louzeiro, Ronny Bergmann, Roland Herzog},
	pdfkeywords={convex analysis, Fenchel conjugate function, Riemannian manifold, Hadamard manifold}
}

\title{Fenchel Duality and a Separation Theorem on Hadamard Manifolds\thanks{MSL was supported by a measure which is co-financed by tax revenue based on the budget approved by the members of the Saxon state parliament. Financial support is gratefully acknowledged.}}
\subtitle{}
\shorttitle{Fenchel Duality on Hadamard Manifolds}

\author{Maurício Silva Louzeiro\thanks{Dongguan University of Technology, School of Computer Science and Technology, Dongguan, Guangdong, China (\email{mauriciosilvalouzeiro@gmail.com}, \orcid{0000-0002-4755-3505}).}
\and
Ronny Bergmann\thanks{Norwegian University of Science and Technology, Department of Mathematical Sciences, NO-7041 Trondheim, Norway (\email{ronny.bergmannn@ntnu.no}, \url{https://www.ntnu.edu/employees/ronny.bergmann}, \orcid{0000-0001-8342-7218}).}
\and
Roland Herzog\thanks{Interdisciplinary Center for Scientific Computing, Heidelberg University, 69120 Heidelberg, Germany (\email{roland.herzog@iwr.uni-heidelberg.de}, \url{https://www.tu-chemnitz.de/mathematik/part_dgl/people/herzog}, \orcid{0000-0003-2164-6575}).}}
\shortauthor{M. Silva Louzeiro, R. Bergmann and R. Herzog}

\dedication{}

\IfFileExists{numapde-uncolorizeHyperrefIfChanges.sty}{\RequirePackage{numapde-uncolorizeHyperrefIfChanges}}{}

\begin{document}
\maketitle

\begin{abstract}
In this paper, we introduce a definition of Fenchel conjugate and Fenchel biconjugate on Hadamard manifolds based on the tangent bundle.
Our definition overcomes the inconvenience that the conjugate depends on the choice of a certain point on the manifold, as previous definitions required.
On the other hand, this new definition still possesses properties known to hold in the Euclidean case.
It even yields a broader interpretation of the Fenchel conjugate in the Euclidean case itself.
Most prominently, our definition of the Fenchel conjugate provides a Fenchel-Moreau Theorem for geodesically convex, proper, lower semicontinuous functions.
In addition, this framework allows us to develop a theory of separation of convex sets on Hadamard manifolds, and a strict separation theorem is obtained.\end{abstract}

\begin{keywords}
convex analysis, Fenchel conjugate function, Riemannian manifold, Hadamard manifold\end{keywords}

\begin{AMS}
\href{https://mathscinet.ams.org/msc/msc2010.html?t=49N15}{49N15}, \href{https://mathscinet.ams.org/msc/msc2010.html?t=90C25}{90C25}, \href{https://mathscinet.ams.org/msc/msc2010.html?t=26B25}{26B25}, \href{https://mathscinet.ams.org/msc/msc2010.html?t=49Q99}{49Q99}
\end{AMS}

\section{Introduction}
\label{section:introduction}

A central concept in convex analysis and related optimization algorithms is the notion of Fenchel duality.
On the other hand, separation theorems for convex sets play an important role for the characterization of functions and their Fenchel conjugate.
Among the vast references on these topics, we mention \cite{BauschkeCombettes:2011:1,EkelandTemam:1999:1,Rockafellar:1970:1,Rockafellar:1974:1,Zalinescu:2002:1,Bot:2010:1,Brezis:2011:1}, all of which consider convex analysis and duality in vector spaces.

The topic of optimization on Riemannian manifolds is currently receiving an increasing amount of attention.
We refer the reader to, \eg, \cite{Udriste:1994:1,Bacak:2014:2,AbsilMahonySepulchre:2008:1,Boumal:2020:1} and \cite[Ch.~6]{Rapcsak:1997:1} for background material.
In this context, a theory of duality on Riemannian manifolds has recently emerged, with particular emphasis on non-smooth problems (\cite{BacakBergmannSteidlWeinmann:2016:1,LellmannStrekalovskiyKoetterCremers:2013:1,WeinmannDemaretStorath:2014:1}) and related algorithms (\cite{BergmannPerschSteidl:2016:1,BergmannChanHielscherPerschSteidl:2016:1, FerreiraLouzeiroPrudente:2020:1}).

To the best of our knowledge, there are up to now only two approaches to Fenchel duality on manifolds: one on Hadamard manifolds due to \cite{AhmadiKakavandiAmini:2010:1}, and one on general Riemannian manifolds proposed in \cite{BergmannHerzogSilvaLouzeiroTenbrinckVidalNunez:2021:1}.
In \cite{AhmadiKakavandiAmini:2010:1}, the authors introduced a Fenchel conjugacy-like concept on complete CAT(0)~spaces (usually called Hadamard spaces), using a quasilinearization in terms of distances as the duality product.
For this purpose, a definition of $o$-dual and $o$-bidual was proposed, where $o$ is a point in the Hadamard space.
The authors then show that this concept possesses several properties of the classical Fenchel conjugate on vector spaces, for instance the biconjugation theorem, and a generalization of the subdifferential characterization.

Recently, we developed in \cite{BergmannHerzogSilvaLouzeiroTenbrinckVidalNunez:2021:1} a theory of duality on Riemannian manifolds~$\cM$ by localizing the Fenchel conjugate similar to \cite{Bertsekas:1978:1}.
As it was the case for \cite{AhmadiKakavandiAmini:2010:1}, this concept also requires the choice of a base point~$m$ on the manifold and it uses duality on the tangent space $\tangent{m}$.
Most of the analysis in \cite{BergmannHerzogSilvaLouzeiroTenbrinckVidalNunez:2021:1} is based on properties of this tangent space as a vector space, and we can generalize many properties of the $m$-Fenchel conjugate to Riemannian manifolds.
We finally derived a generalization of the so-called Chambolle-Pock algorithm (\cite{PockCremersBischofChambolle:2009:1,ChambollePock:2011:1}) for the minimization of $f(p) + g(\Lambda p)$, where $f$ is defined on $\cM$, $g$ is defined on another Riemannian manifold~$\cN$, and $\Lambda \colon \cM \to \cN$.
The algorithm generalizes a concept from \cite{Valkonen:2014:1} and employs a linearization of $\Lambda$ at a point $m$ as well as the $n$-Fenchel conjugate for a second base point $n \in \cN$.
The convergence results rely on the convexity of the pull-back of $g$ onto the tangent space, \ie, on the convexity of the composition $g \circ \exponential{m}$.

In this paper we introduce a competing definition of a Fenchel conjugate on Hadamard manifolds.
Our new definition differs from the one in \cite{BergmannHerzogSilvaLouzeiroTenbrinckVidalNunez:2021:1} in two important ways.
First, the conjugate of a function $F \colon \cM \to \eR$ is defined on the entire cotangent bundle, not just on the cotangent space at a particular base point.
Second, we do not pull $F$ back to the tangent space.
We also define the Fenchel biconjugate, which is---similar as in the competing approaches---again a function defined on the manifold.

Our new concept of duality possesses similar properties as those proved for the $m$-Fenchel conjugate in \cite{BergmannHerzogSilvaLouzeiroTenbrinckVidalNunez:2021:1}.
These include, in particular, the characterization of the subdifferential in terms of the Fenchel conjugate as well as the biconjucation theorem, also known as Fenchel--Moreau Theorem.
The main difference is that these results hold under more natural assumptions, notably geodesic convexity of the function under consideration, rather than the convexity of its pull-back to the tangent space.
We thus envision our work to give new insight into potentially new algorithmic concepts for geodesically convex optimization problems on Riemannian manifolds, which have numerous recent applications, \eg, in signal and image processing \cite{FletcherJoshi:2007:1,BergmannFitschenPerschSteidl:2018:1,BergmannGousenbourger:2018:1} and statistics and machine learning \cite{JourneeBachAbsilSepulchre:2010:1,Wiesel:2012:1,Vandereycken:2013:1,HosseiniSra:2015:1,AllenZhuGargLiOliveiraWigderson:2018:1,GoyalShetty:2019:1}.

We would like to emphasize that our definition of Fenchel conjugate provides a broader understanding of the concept even for functions~$f$ defined on a vector space~$V$.
While classically, the conjugate $f^*$ is a function defined on $V^*$, we obtain here a conjugate $f^*$ defined on the cotangent bundle~$V \times V^*$, for which $f^*(0,\cdot)$ agrees with the classical definition and any section $f^*(x,\cdot)$ determines all other sections $f^*(y,\cdot)$; see \cref{rem:eucli.case}.
In fact, it was one of the main conceptual difficulties to recognize that in the case of manifolds, the conjugate function should contain information not only from a single space of cotangent directions but rather from all cotangent spaces in order to recover the Biconjugation \cref{theo:main} under the same, natural assumptions as in the case of vector spaces.

An additional result in this paper is a theorem regarding the separation of convex sets on Hadamard manifolds by geodesic hyperplanes in the cotangent bundle.
This generalizes a well known separation theorem from vector spaces to Hadamard spaces.

The remainder of the paper is organized as follows.
In \cref{sec:Preliminaries} we recall a number of classical results from convex analysis in Hilbert spaces.
In an effort to make the paper self-contained, we also briefly state the required concepts from differential geometry and convex analysis on Hadamard manifolds.
\Cref{section:fenchel.conjugate} is devoted to the development of the new notion of Fenchel conjugation for functions defined on Hadamard manifolds.
Leveraging the concept, we extend some classical results from convex analysis to manifolds, like the Fenchel--Moreau Theorem (also known as the Biconjugation Theorem) and the characterization of the subdifferential in terms of the conjugate function.
In \cref{section:separation-theory}, we introduce a theory of separation of convex sets on Hadamard manifolds, which leads to a strict separation theorem.
Finally, we give some conclusions and further remarks on future research in \cref{sec:Conclusions}.

\section{Preliminaries on Convex Analysis and Differential Geometry}
\label{sec:Preliminaries}

In this section we review some well known results from convex analysis in Hilbert spaces, which serve as the standard for comparison for the new results to be developed in \cref{section:fenchel.conjugate}.
We emphasize that the results collected here are valid in more general contexts, but we do not strive for full generality.
We also revisit necessary concepts from differential geometry as well as the intersection of both topics, convex analysis on Riemannian manifolds, including its subdifferential calculus.

Throughout this paper we denote the extended line as $\eR \coloneqq \R \cup \{\pm \infty\}$.
We shall use the usual convention $-(-\infty) = +\infty$ and $-(+\infty) = -\infty$.

\subsection{Convex Analysis}
\label{subsec:Convex_Analysis}

In this subsection let~$\cX$ be a Hilbert space with inner product $\inner{\cdot}{\cdot}\colon\cX \times \cX \to \eR$ and duality pairing $\dual{\cdot}{\cdot}\colon \cX^*\times \cX \to \eR$.
Here, $\cX^*$ denotes the dual space of $\cX$.
For standard definitions like \emph{closedness, properness, lower semicontinuity (lsc)} and \emph{convexity} of a function~$f\colon \cX \to \eR$, we refer the reader, \eg, to the textbooks \cite{Rockafellar:1970:1,BauschkeCombettes:2011:1}.

\begin{definition}
	\label{def:Classical_Fenchel_conjugate}%
	The \emph{Fenchel conjugate} of a function~$f\colon \cX \to\eR$ is defined as the function $f^*\colon\cX^* \to \eR$ such that
	\begin{equation*}
		\label{eq:Classical_Fenchel_conjugate}
		f^*(x^*)
		\coloneqq
		\sup_{x \in \cX} \paren[big]\{\}{\dual{x^*}{x} - f(x)}.
	\end{equation*}
\end{definition}
We mention that some authors define equivalently $f^* \colon \cX \to \eR$, replacing the duality pairing $\dual{x^*}{x}$ by the inner product $\inner{x^*}{x}$.
A similar statement applies to the definition of the subdifferential in \cref{def:Classiscal_convex_subdifferential} below.

We recall some properties of the Fenchel conjugate function in Hilbert spaces in the following lemma.

\begin{lemma}[{\cite[Ch.~13]{BauschkeCombettes:2011:1}}]
	\label{lemma:Classical_properties_Fenchel_conjugate}%
	Let~$f,g\colon \cX \to\eR$ be proper functions, $\alpha \in \R$, $\lambda > 0$ and $z \in \cX$.
	Then the following statements hold.
	\begin{enumerate}
		\item\label{item:Classical_properties_conjugate}
			$f^*$ is convex and lsc.
		\item\label{item:Classical_properties_inequality}
			If $f(x) \le g(x)$ for all $x \in \cX$, then $f^*(x^*) \ge g^*(x^*)$ for all $x^* \in \cX^*$.
		\item\label{item:Classical_properties_adding_alpha}
			If $g(x) = f(x)+\alpha$ for all $x \in \cX$, then $g^*(x^*) = f^*(x^*) - \alpha$ for all $x^* \in \cX^*$.
		\item\label{item:Classical_properties_lambda}
			If $g(x) = \lambda f(x)$ for all $x \in \cX$, then $g^*(x^*) = \lambda f^*(x^*/\lambda)$ for all $x^* \in \cX^*$.
		\item\label{item:Classical_properties_Fenchel_shift}
			If $g(x) = f(x+z)$ for all $x \in \cX$, then $g^*(x^*) = f^*(x^*) - \dual{x^*}{z}$ for all $x^* \in \cX^*$.
		\item\label{item:Classical_properties_Fenchel-Young}
			The \emph{Fenchel--Young inequality} holds, \ie, for all $(x,x^*) \in \cX \times \cX^*$ we have
			\begin{equation*}
				\label{eq:Classical_Fenchel_Young_inequality}
				\dual{x^*}{x} \le f(x) + f^*(x^*).
			\end{equation*}
	\end{enumerate}
\end{lemma}

We now recall some results related to the definition of the subdifferential of a proper function.

\begin{definition}[{\cite[Def.~16.1]{BauschkeCombettes:2011:1}}]
	\label{def:Classiscal_convex_subdifferential}%
	Let $f\colon \cX \to\eR$ be a proper function.
	Its subdifferential is defined as
	\begin{equation*}
		\partial f(x)
		\coloneqq
		\setDef[auto]{x^* \in \cX^*}{f(z) \ge f(x) + \dual{x^*}{z - x} \quad \text{for all } z \in \cX}.
	\end{equation*}
\end{definition}

\begin{theorem}[{\cite[Prop.~16.9]{BauschkeCombettes:2011:1}}]
	\label{thm:Classical_characterization_subdifferential}%
	Let $f\colon \cX \to \eR$ be a proper function and~$x \in \cX$.
	Then $x^* \in \partial f(x)$ holds if and only if
	\begin{equation*}
		f(x)+f^*(x^*) = \dual{x^*}{x}.
	\end{equation*}
\end{theorem}

The Fenchel biconjugate~$f^{**} \colon \cX \to\eR$ of a function~$f \colon \cX \to \eR$ is given by
\begin{equation}
	\label{eq:Classical_biconjugate_function}
	f^{**}(x)
	=
	(f^*)^*(x)
	=
	\sup_{x^* \in \cX^*} \paren[big]\{\}{\dual{x^*}{x} -f^*(x^*)}.
\end{equation}
It satisfies $f^{**}(x) \le f(x)$ for all $x \in \cX$; see for instance \cite[Prop.~13.14]{BauschkeCombettes:2011:1}.

We conclude this section with the following result known as the Fenchel--Moreau or Biconjugation Theorem.

\begin{theorem}[{\cite[Thm.~13.32]{BauschkeCombettes:2011:1}}]
	\label{thm:Classical_Fenchel_Moreaou_theorem}%
	Given a proper function~$f \colon \cX \to \eR$, the equality~$f^{**}(x) = f(x)$ holds for all $x \in \cX$ if and only if $f$ is lsc and convex.
	In this case $f^*$ is proper as well.
\end{theorem}

\subsection{Differential Geometry on Riemannian Manifolds}
\label{subsec:Differential_Geometry}

This section is devoted to the collection of necessary concepts from differential geometry.
For details concerning the subsequent definitions, the reader may wish to consult \cite{DoCarmo:1992:1,Lee:2003:1,Jost:2017:1}.

Suppose that~$\cM$ is an~$n$-dimensional connected, smooth manifold.
The tangent space at $p \in \cM$ is a vector space of dimension~$n$ and it is denoted by $\tangent{p}$.
Its dual space is denoted by~$\cotangent{p}$ and it is called the \emph{cotangent space} to $\cM$ at $p$.
The duality product between $X \in \tangent{p}$ and $\xi \in \cotangent{p}$ is denoted by $\dual{\xi}{X} = \xi(X) \in \R$.

The disjoint union of all tangent respectively cotangent spaces, \ie,
\begin{equation*}
	\tangentBundle
	\coloneqq
	\bigcup_{p \in \cM} \paren\{\}{p} \times \tangent{p}
	\quad
	\text{and}
	\quad
	\cotangentBundle
	\coloneqq
	\bigcup_{p \in \cM} \paren\{\}{p} \times \cotangent{p}
\end{equation*}
is called the \emph{tangent bundle} respectively the \emph{cotangent bundle} of~$\cM$.
Both are smooth manifolds of dimension~$2n$.

We suppose that $\cM$ is equipped with a Riemannian metric, \ie, a smoothly varying family of inner products on the tangent spaces $\tangent{p}$.
The metric at $p \in \cM$ is denoted by $\riemannian{\cdot}{\cdot}[p] \colon \tangent{p}\times \tangent{p} \to\R$ and we write $\riemanniannorm{\cdot}[p]$ for the associated norm in $\tangent{p}$.
For simplicity we shall omit the index $p$ when no ambiguity arises.
The Riemannian metric furnishes a linear bijective correspondence between the tangent and cotangent spaces via the Riesz map and its inverse, the so-called \emph{musical isomorphisms}; see \cite[Ch.~8]{Lee:2003:1}.
They are defined as
\begin{alignat}{2}
	\label{eq:Flat_isomorphism}
	\flat\colon\tangent{p}\ni X
	&
	\mapsto X^\flat \in \cotangent{p},
	&
	\quad
	\dual{X^\flat}{Y}
	&
	=
	\riemannian{X}{Y}[p],
	\quad
	\text{for all } Y \in \tangent{p}
	,
	\intertext{and its inverse,}
	\label{eq:Sharp_isomorphism}
	\sharp\colon\cotangent{p}\ni\xi
	&
	\mapsto\xi^\sharp \in \tangent{p},
	&
	\quad
	\riemannian{\xi^\sharp}{Y}[p]
	&
	=
	\dual{\xi}{Y},
	\quad
	\text{for all } Y \in \tangent{p}
	.
\end{alignat}
The $\sharp$-isomorphism further introduces an inner product and an associated norm on the cotangent space $\cotangent{p}$, which we will also denote by~$\riemannian{\cdot}{\cdot}[p]$ and $\riemanniannorm{\cdot}[p]$, since it is clear which inner product or norm we refer to based on the respective arguments.

The tangent vector of a curve~$\gamma \colon I \to \cM$ defined on some open interval~$I \subseteq \R$ is denoted by $\dot \gamma(t)$.
A curve is said to be geodesic if $\nabla_{\dot \gamma(t)} \dot \gamma(t) = 0$ holds for all $t \in I$, where $\nabla$ denotes the Levi-Cevita connection, cf.\ \cite[Ch.~2]{DoCarmo:1992:1} or \cite[Thm.~4.24]{Lee:2018:1}.
As a consequence, geodesic curves have constant speed.
We say that a geodesic $\gamma \colon [0,1] \subset \R \to \cM$ connects $p$ to $q$ if $\gamma(0) = p$ and $\gamma(1) = q$ holds.
Note that a geodesic connecting $p$ to $q$ need not always exist, and if it exists, it need not be unique.
Geodesics might also be of different lengths. 
If a unique shortest geodesic connecting $p$ and $q$ exists, we denote it by $\geodesic<a>{p}{q}$.
Moreover, given $(p,X) \in \tangentBundle$, we denote by~$\geodesic{p}{X} \colon I \to \cM$, with $I \subseteq \R$ being a suitable open interval containing~$0$, the geodesic starting at~$p$ with~$\dot{\gamma}_{p,X}(0) = X$.
We denote the subset of $\tangent{p}$ for which these geodesics are well defined until $t = 1$ by $\cG_p$.
Recall that a Riemannian manifold~$\cM$ is said to be (geodesically) \emph{complete} if $\cG_p = \tangent{p}$ holds for some, and equivalently for all~$p \in \cM$.

The Riemannian distance between $p$ and $q$ in $\cM$, defined as the infimum of the length over all piecewise smooth curve segments from $p$ to $q$, is denoted by $d(p,q)$.
The metric topology it induces agrees with the original topology on $\cM$.
By the Hopf-Rinow theorem, $\cM$ is geodesically complete if and only if it complete in the sense of metric spaces.

The \emph{exponential map} is defined as the function~$\exponential{p} \colon \cG_p\to\cM$ with~$\exponential{p} X \coloneqq \geodesic{p}{X}(1)$.
Note that~$\exponential{p}(tX) = \geodesic{p}{X}(t)$ holds for every $t \in [0,1]$.
We further introduce the set $\cG'_p \subseteq \tangent{p}$ as some open set such that $\exponential{p} \colon \cG'_p \to \exponential{p}(\cG'_p) \subseteq \cM$ is a diffeomorphism.
The \emph{logarithmic map} is defined as the inverse of the exponential map, \ie, $\logarithm{p} \colon \exponential{p}(\cG'_p) \to \cG'_p \subseteq \tangent{p}$.

In the particular case of a \emph{Hadamard manifold}, \ie, a manifold which is simply connected and complete and whose  sectional curvature is nonpositive everywhere, the geodesics connecting any two distinct points exist and are unique; see \cite[p.~10]{Bacak:2014:2}.
In this case, the exponential and logarithmic maps are defined globally, \ie,~ $\cG_p = \tangent{p}$ holds for all $p \in \cM$.
Moreover, the distance function “is at least as convex as in the Euclidean plane” \cite[p.6]{Bacak:2014:2} and the squared distance function $d(\cdot,p)^2$ is even strongly convex \cite[Rem.~2.2.2]{Bacak:2014:2}.
These properties of Hadamard manifolds make these spaces particularly amenable for the study of convexity properties.

\subsection{Convex Analysis on Hadamard Manifolds}
\label{subsec:Convex_analysis_on_Riemannian_manifolds}

Throughout this subsection, $\cM$ is assumed to be a Hadamard manifold and we recall the basic concepts of convex analysis on $\cM$.
The central idea is to replace straight lines in the definition of convex sets in Hilbert spaces by geodesics.

\begin{definition}[{\cite[Def.~IV.5.9, Def.~IV.5.1]{Sakai:1996:1}}] \hfill
	\label{def:F_geodesically_convex}%
	\begin{enumerate}
		\item
			A function~$F\colon \cM\to\eR$ is \emph{proper} if $\dom F \coloneqq \setDef{p \in \cM}{F(p) < \infty} \neq \emptyset$ and $F(p) > -\infty$ holds for all~$p \in \cM$.
		\item
			\label{item:F_geodesically_convex}
			A function~$F\colon\cM\to \eR$ is \emph{convex} if, for all~$p,q \in \cM$, the composition $F \circ \geodesic<a>{p}{q} \colon [0,1] \subset \R \to \eR$ is a convex function on $[0,1]$ in the classical sense.
		\item
			\label{item:epigraph}
			The \emph{epigraph} of a function $F \colon \cM \to \eR$ is defined as
			\begin{equation}
				\label{eq:Epigraph}
				\epi F
				\coloneqq
				\setDef{(p,\alpha) \in \cM \times \R}{F(p) \le \alpha}.
			\end{equation}
		\item
			\label{item:lsc}
			A \emph{proper} function~$F\colon \cM \to\eR$ is called \emph{lower semicontinuous (lsc)} if $\epi F$ is closed.
		\item
			\label{item:conv.subs}
			A subset $\cC \subseteq \cM$ is said to be \emph{convex} if for any two points $p, q \in \cC$, the unique geodesic of $\cM$ connecting $p$
			to $q$ lies completely in $\cC$.
	\end{enumerate}
\end{definition}

An interesting observation here is, that geodesic balls around $p\in \cM$ of radius $r\geq 0$, i.e.
\begin{equation*}
\cB_p(r)
\coloneqq
\setDef{q\in \cM}{\text{ there exists } X \in \tangent{p} \text{ with } \norm{X}_p \leq r \text{ such that } q = \exponential{p} X}
\end{equation*}
are convex sets.

We now recall the notion of the subdifferential of a geodesically convex function.
\begin{definition}[{\cite{FerreiraOliveira:1998:1}, \cite[Def.~3.4.4]{Udriste:1994:1}}]
	\label{def:Subdifferential}%
	The \emph{subdifferential}~$\partial F$ at a point~$p \in \cM$ of a proper, convex function~$F\colon\cM\to\eR$ is given by
	\begin{equation*}\label{eq:Subdifferential_manifold}
		\partial F(p)
		\coloneqq
		\setDef[big]{\xi \in \cotangent{p}}{F(q) \ge F(p) + \dual{\xi}{\logarithm{p} q} \quad \text{for all } q \in \cM}
		.
	\end{equation*}
\end{definition}
As was mentioned already for the Hilbert space case, the subdifferential is sometimes defined equivalently as a subset of the tangent space, and the duality pairing $\dual{\xi}{\logarithm{p} q}$ is replaced by an inner product.

When $\cC \subseteq \cM$ is nonempty, convex and closed, it was proved in \cite{FerreiraOliveira:2002:1} that for each point $p \in \cM$, there is a unique point $\hat p \in \cC$ satisfying $d(p,\hat p) \le d(p, q)$ for all $q \in \cC$.
In this case, $\hat p$ is called the \emph{projection} of $p$ onto $\cC$ and it is denoted by $\proj{\cC}(p)$.
We require the following result from \cite[Cor.~3.1]{FerreiraOliveira:2002:1}.

\begin{theorem}\label{theo:ineq.proj}
	Suppose that $\cC \subseteq \cM$ a nonempty, convex and closed set and $p \in \cM$.
	Then the following inequality holds,
	\begin{equation*}
		\riemannian[big]{\logarithm{\proj{\cC}(p)} p}{\logarithm{\proj{\cC}(p)} q}
		\le
		0
		\quad \text{for all } q \in \cC
		.
	\end{equation*}
\end{theorem}
\begin{corollary}\label{coro:ineq.proj}
	Let $F \colon \cM \to \eR$ be a proper lsc convex function and $(p,s) \not \in \epi F$.
	Then the projection $\proj{\epi F}(p,s) \eqqcolon (\hat p,\hat s)$ exists and the following inequality holds,
	\begin{equation*}
		\riemannian[big]{\logarithm{\hat p} p}{\logarithm{\hat p} q} + (s-\hat s)(r-\hat s)
		\le
		0
		\quad \text{for all } (q,r) \in \epi F
		.
	\end{equation*}
\end{corollary}
\begin{proof}
	Since $F \colon \cM \to \eR$ is a proper lsc convex function, $\epi F \subset \cM \times \R$ is a nonempty closed convex set, where $\cM \times \R$ is equipped with the product metric.
	Hence, using~\cref{theo:ineq.proj} with $\cC = \epi F$ we get the desired inequality.
\end{proof}

A geodesic triangle $\geodesicTriangle{p_0}{p_1}{p_2}$ of a Hadamard manifold is the set consisting of three distinct points $p_0,p_1,p_2$ called the vertices and three geodesics $\geodesic<a>{p_0}{p_1}$, $\geodesic<a>{p_1}{p_2}$, $\geodesic<a>{p_2}{p_0}$.
The proof of the following theorem can be found in \cite[Thm.~2.2]{FerreiraOliveira:2002:1}.
\begin{theorem}\label{theo:cosine.law}
	Suppose that $\geodesicTriangle{p_0}{p_1}{p_2}$ a geodesic triangle.
	Then,
	\begin{align}
		d^2(p_{i \ominus 1},p_i) - 2 \riemannian[big]{\logarithm{p_i} p_{i \ominus 1}}{\logarithm{p_i} p_{i \oplus 1}} + d^2(p_{i \oplus 1},p_i)
		&
		\le
		d^2(p_{i \ominus 1},p_{i \oplus 1})
		,
		\label{eq2:theo:cosine.law}
		\\
		\riemannian[big]{\logarithm{p_i} p_{i \ominus 1}}{\logarithm{p_i} p_{i \oplus 1}} + \riemannian[big]{\logarithm{p_{i \oplus 1}} p_{i \ominus 1}}{\logarithm{p_{i \oplus 1}} p_{i}}
		&
		\ge
		d^2(p_{i \oplus 1}, p_i)
		,
		\label{eq3:theo:cosine.law}
	\end{align}
	for $i = 0,1,2$, where the indices $i \ominus 1$ and $i \oplus 1$ are meant modulo~$3$.
\end{theorem}

\section{Fenchel Conjugate on Hadamard Manifolds}
\label{section:fenchel.conjugate}

In this section we introduce new definitions of the Fenchel conjugate and Fenchel biconjugate for extended real-valued functions defined on Hadamard manifolds.
Using these definitions, we can extend fundamental properties from the Euclidean to the Riemannian setting.
In \cref{subsec:application}, we elaborate on potential applications and provide a concrete example of the Fenchel conjugate of a function on the manifold of symmetric, positive definite matrices.

\subsection{Fenchel Conjugate}
\label{subsec:Fenchel-Conjugate}

We begin with a new definition of the Fenchel conjugate function on Hadamard manifolds.
\begin{definition}\label{def:Fenconj}
	Let $F \colon \cM \to \eR$.
	The Fenchel conjugate of $F$ is the function $F^* \colon \cotangentBundle \to \eR$ defined by
	\begin{equation}\label{eq:conjugate}
		F^*(p,\xi)
		\coloneqq
		\sup_{q \in \cM} \paren[big]\{\}{\dual{\xi}{\logarithm{p} q} - F(q)}
		\quad \text{for } (p,\xi) \in \cotangentBundle
		.
	\end{equation}
\end{definition}
As was mentioned in the introduction, this definition differs from \cite[Def.~3.1]{BergmannHerzogSilvaLouzeiroTenbrinckVidalNunez:2021:1} in two important ways.
First, $F^*$ is defined on the entire cotangent bundle, not just on the cotangent space at a particular base point.
Second, we do not pull $F$ back to the tangent space.

\begin{remark}\label{rem:conjugate.dom}
	Suppose that $F \colon \cM \to \eR$ is a proper function.
	Since $\dual{\xi}{\logarithm{p} q} - F(q) = -\infty$ holds for all $q \not \in \dom F$, we have
	\begin{equation*}
		F^*(p,\xi)
		=
		\sup_{q \in \dom F} \paren[big]\{\}{\dual{\xi}{\logarithm{p} q} - F(q)}
		\quad \text{for all } (p,\xi) \in \cotangentBundle
		.
	\end{equation*}
\end{remark}

\begin{remark}\label{rem:previous-def}
	For each $(p,\xi) \in \cotangentBundle$, $F^*(p,\xi)$ agrees with the $p$-Fenchel conjugate $F^*_p(\xi)$ introduced in \cite[Def.~3.1]{BergmannHerzogSilvaLouzeiroTenbrinckVidalNunez:2021:1}.
	For convenience, let us recall that $F^*_p \colon \cotangent{p} \to \eR$ is defined as
	\begin{equation}\label{eq:old-def}
		F^*_p(\xi)
		=
		\sup_{X \in \tangent{p}[\cM]} \paren[big]\{\}{\dual{\xi}{X} - F(\exponential{p} X)}
		\quad \text{for } \xi \in \cotangent{p}[\cM]
		.
	\end{equation}
	The equality $F^*(p,\xi) = F^*_p(\xi)$ follows immediately from the relation $X = \logarithm{p} q \Leftrightarrow q = \exponential{p} X$ on Hadamard manifolds.
\end{remark}
\begin{remark}\label{rem:eucli.case}
	We also observe that in case $\cM$ is the Euclidean space $\R^n$, \cref{def:Fenconj} becomes
	\begin{align}
		\label{eq:redundancy_in_the_conjugate_on_Rn}
		F^*(p,\xi)
		&
		=
		\sup_{q \in \R^n} \paren[big]\{\}{\dual{\xi}{q} - F(q)} - \dual{\xi}{p}
		\notag
		\\
		&
		=
		F^*(\xi) - \dual{\xi}{p}
		\quad \text{for all } (p,\xi) \in \R^n \times \R^n
		.
	\end{align}
	Due to \cref{lemma:Classical_properties_Fenchel_conjugate}~\cref{item:Classical_properties_Fenchel_shift} this is the Fenchel conjugate of $F(p+\xi)$.
	Similar but not identical to $F_p^*$, cf.~\eqref{eq:old-def} or \cite[Def.~3.1]{BergmannHerzogSilvaLouzeiroTenbrinckVidalNunez:2021:1},
	we can recover the classical (Euclidean) case.
	The domain of the Fenchel conjugate here is larger than in the Euclidean case.
	If we set $p = 0$ in the first argument, \ie consider the function
	$F^*(0,\cdot)$ only in its second argument, then we obtain $F^*(0,\cdot) = F_0^* = F^*$.
\end{remark}

\begin{example}
	\label{example:distance}
	Let $p \in \cM$ be arbitrary but fixed and let $F \colon \cM \to \R$ be defined by $F(q) = d(p,q)$.
	Due to the fact that $d(p,q) = \norm{\logarithm{p} q}$ holds, we obtain from \cref{def:Fenconj} the following representation of $F^*$:
	\begin{equation}\label{eq:conjugate.distance}
		F^*(p,\xi)
		=
		\sup_{q \in \cM} \paren[big]\{\}{\dual{\xi}{\logarithm{p} q} - \norm{\logarithm{p} q}}
		\quad \text{for } (p,\xi) \in \cotangentBundle
		.
	\end{equation}
	For every $\xi \in \cotangent{p}[\cM]$ with $\norm{\xi} \le 1$, the following inequalities hold:
	\begin{align*}
		0
		=
		\dual{\xi}{\logarithm{p} p} - \norm{\logarithm{p} p}
		&
		\le
		\sup_{q \in \cM} \paren[big]\{\}{\dual{\xi}{\logarithm{p} q} - \norm{\logarithm{p} q}}
		\\
		&
		\le
		\sup_{q \in \cM} \paren[big]\{\}{\norm{\xi} \norm{\logarithm{p} q} - \norm{\logarithm{p} q}}
		\le
		0
		.
	\end{align*}
	Hence, \eqref{eq:conjugate.distance} implies $F^*(p,\xi) = 0$ whenever $\norm{\xi} \le 1$.
	On the other hand, if $\norm{\xi} > 1$ holds, then
	\begin{align*}
		F^*(p,\xi)
		&
		=
		\sup_{q \in \cM} \paren[big]\{\}{\dual{\xi}{\logarithm{p} q} - \norm{\logarithm{p} q}}
		\\
		&
		=
		\sup_{X \in \tangent{p}} \paren[big]\{\}{\dual{\xi}{X} - \norm{X} }
		\ge
		\sup_{\lambda > 0} \paren[big]\{\}{\dual{\xi}{\lambda \, \xi^\sharp} - \norm{\lambda \, \xi^\sharp}}
		\\
		&
		=
		\sup_{\lambda > 0} \paren[big]\{\}{\lambda \, (\norm{\xi}^2 - \norm{\xi})}
		=
		+ \infty
		.
	\end{align*}
	Therefore, the Fenchel conjugate $F^*\colon \cotangentBundle \to \eR$ of $F = d(p,\cdot)$ is given by
	\begin{equation*}
		F^*(p,\xi)
		=
		\begin{cases}
			0 & \text{if } \norm{\xi} \le 1,
			\\
			+\infty & \text{if } \norm{\xi} > 1 .
		\end{cases}
	\end{equation*}
\end{example}

\begin{example}
	\label{example:distance_squared}
	Let $p \in \cM$ be arbitrary but fixed and let $F \colon \cM \to \R$ be defined by $F(q) = \frac{1}{2}d^2(p,q)$.
	Then we have
	\begin{align*}
		F^*(p,\xi)
		&
		=
		\sup_{q \in \cM} \paren[Big]\{\}{\dual{\xi}{\logarithm{p} q} - \frac{1}{2}\riemanniannorm{\logarithm{p} q}^2}
		\\
		&
		=
		\sup_{X \in \tangent{p}} \paren[Big]\{\}{\dual{\xi}{X} - \frac{1}{2}\riemanniannorm{X}^2}
		=
		\frac{1}{2}\riemanniannorm{\xi}^2
		.
	\end{align*}
	In particular, we obtain $F^*(p,\xi) = F^*(p,-\xi)$ for all $(p,\xi) \in \cotangentBundle$.
	In addition, the following property holds:
	\begin{equation*}
		F^*\paren[big](){p,[\logarithm{p} q]^{\flat}}
		=
		\frac{1}{2}\riemanniannorm{[\logarithm{p} q]^{\flat}}^2
		=
		\frac{1}{2}d^2(p,q)
		=
		F(q)
		\quad \text{for all } q \in \cM
		.
	\end{equation*}
\end{example}

In comparison with the classical conjugate on $\R^n$, it appears unusual that $F^*$ from \cref{def:Fenconj} depends on two arguments, $p$ and $\xi$.
One might expect there to be some redundancy in the definition.
Indeed, this redundancy has already been observed in \eqref{eq:redundancy_in_the_conjugate_on_Rn} for the Euclidean setting.
We now explore it in the Riemannian case.

To this end, we consider the following equivalence relation on the cotangent bundle $\cotangentBundle$:
\begin{equation}
	\label{eq:equivalence_relation_on_cotangent_bundle}
	(p,\xi) \sim (p',\xi')
	\quad
	\text{if and only if}
	\quad
	\dual{\xi}{\logarithm{p} q} = \dual{\xi'}{\logarithm{p'} q}
	\text{ holds for all }
	q \in \cM
	.
\end{equation}
The equivalence class of $(p,\xi) \in \cotangentBundle$, denoted by $[(p,\xi)]$, is
\begin{equation}\label{eq:equivalence-class}
	[(p,\xi)]
	=
	\setDef{(p',\xi') \in \cotangentBundle}{\dual{\xi}{\logarithm{p} q} = \dual{\xi'}{\logarithm{p'} q} \text{ for all } q \in \cM}
	.
\end{equation}
Note that $F^*(p',\xi') = F^*(p,\xi)$ holds for all $(p',\xi') \in [(p,\xi)]$.
Observe as well that when $\cM$ is the Euclidean space $\R^n$, the equivalence class of $(p,\xi) \in \cotangentBundle$ becomes
\begin{align*}
	[(p,\xi)]
	&
	=
	\setDef{(p',\xi') \in \R^n \times \R^n}{\dual{\xi}{q - p} = \dual{\xi'}{q - p'} \text{ for all } q \in \R^n}
	\\
	&
	=
	\setDef{(p',\xi') \in \R^n \times \R^n}{\dual{\xi-\xi'}{q} = \dual{\xi}{p} - \dual{\xi'}{p'} \text{ for all } q \in \R^n}
	\\
	&
	=
	\setDef{(p',\xi') \in \R^n \times \R^n}{\xi' = \xi, \; \dual{\xi'}{p'} = \dual{\xi}{p}}
	\\
	&
	=
	\setDef{(p',\xi) \in \R^n \times \R^n}{\dual{\xi}{p'} = \dual{\xi}{p}}
	,
\end{align*}
which describes a hyperplane in $\R^n$ when $\xi \neq 0$.
The following example illustrates that the equivalence class defined in \eqref{eq:equivalence-class} is, in general, not a singleton even in non-Euclidean manifolds.

\begin{example}\label{ex:nosingleton}
	We denote by $\cM = \spd{n}$ the cone of real, symmetric positive definite matrices of size $n \times n$.
	Its tangent space (at any point) can be identified with $\symmetric{n}$, the space of symmetric $n \times n$-matrices.
	The manifold $\cM$ is endowed with the affine invariant Riemannian metric, which at $A \in \spd{n}$ is given by
	\begin{equation}\label{eq:spd_metric}
		\riemannian{X}{Y}[A]
		\coloneqq
		\trace (X A^{-1} Y A^{-1})
		\quad
		\text{for }
		X, Y \in \tangent{A}[\cM]
		.
	\end{equation}
	$\cM$ is a Hadamard manifold; see for instance \cite[Ch.~XII, Thm.~1.2, p.~325]{Lang:1999:1}.
	When we identify the cotangent space with $\symmetric{n}$ via the duality $\dual{\xi}{Y} \coloneqq \trace(\xi \, Y)$, then the 'flat' isomorphism~$\flat$ at $A$ is given by
	\begin{equation}
		\label{eq:spd_flat_isomorphism}
		X^\flat
		=
		A^{-1} X A^{-1}
	\end{equation}
	since $\dual{X^\flat}{Y} = \trace(X^\flat Y) = \trace(A^{-1} X A^{-1} Y) = \trace(X A^{-1} Y A^{-1}) = \riemannian{X}{Y}[A]$ holds for all $Y \in \symmetric{n}$.

	The logarithmic map $\logarithm{A} \colon \cM \to \tangent{A}[\cM]$ is given by
	\begin{equation} \label{eq:spd_logarithmic_map}
		\logarithm{A} B
		=
		A^{1/2} \matrixLogarithm\paren[big](){A^{-1/2} B \, A^{-1/2}} \, A^{1/2}
		\quad
		\text{for }	A, B \in \cM
		,
	\end{equation}
	where $\cdot^{1/2}$ and $\matrixLogarithm$ denote the matrix square root and matrix logarithm of symmetric positive definite matrices, respectively.
	We refer the reader, for instance, to \cite[Thms.~1.29 and 1.31]{Higham:2008:1}.

	Suppose that $A \in \spd{n}$ is arbitrary but fixed and consider the particular cotangent vector $A^\flat = A^{-1} \in \symmetric{n}$.
	Using \eqref{eq:spd_logarithmic_map}, we evaluate
	\begin{align*}
		\dual[big]{A^\flat}{\logarithm{A} B}
		&
		=
		\trace \paren[big](){A^{-1} A^{1/2} \matrixLogarithm\paren[big](){A^{-1/2} B \, A^{-1/2}} \, A^{1/2}}
		\\
		&
		=
		\trace \matrixLogarithm\paren[big](){A^{-1/2} B \, A^{-1/2}}
		\\
		&
		=
		\trace \matrixLogarithm A^{-1/2}
		+
		\trace \matrixLogarithm B
		+
		\trace \matrixLogarithm A^{-1/2}
		\\
		&
		=
		\trace \matrixLogarithm B
		-
		\trace \matrixLogarithm A
	\end{align*}
	for any $B \in \spd{n}$.
	Here we also used that $\trace \matrixLogarithm (CD) = \trace \matrixLogarithm C + \trace \matrixLogarithm D$ holds for positive definite matrices $C$ and $D$ as well as $\matrixLogarithm C^{-1} = - \matrixLogarithm C$.

	Now choose any orthogonal matrix $V$ and set $\widehat A \coloneqq V \, A \, V^{-1}$.
	Then the same reasoning as above shows
	\begin{equation*}
		\dual[big]{\widehat A^\flat}{\logarithm{\widehat A} B}
		=
		\trace \matrixLogarithm B
		-
		\trace \matrixLogarithm \widehat A
		=
		\trace \matrixLogarithm B
		-
		\trace \matrixLogarithm A
		.
	\end{equation*}
	The second equality follows from the fact that $A$ and $\widehat A$ have the same eigenvalues and thus the same is true for $\matrixLogarithm A$ and $\matrixLogarithm \widehat A$.

	We conclude that $(A,A^\flat)$ and $(\widehat A,\widehat A^\flat)$ belong to the same equivalence class \wrt the relation \eqref{eq:equivalence_relation_on_cotangent_bundle}.
	Therefore, the equivalence classes \eqref{eq:equivalence-class} are not, in general, singletons.
\end{example}

The following results establishes a property which relates the Fenchel conjugate evaluated in elements of $\cotangentBundle$ whose base points are not necessarily the same.
\begin{proposition}
	Let $F \colon \cM \to \overline{\R}$ and $p,p' \in \cM$.
	Then the following inequality holds:
	\begin{equation*}
		F^*\paren[big](){p, [\logarithm{p} p']^{\flat}}
		\ge
		F^*\paren[big](){p', [-\logarithm{p'} p]^{\flat}} + d^2(p,p')
		.
	\end{equation*}
\end{proposition}
\begin{proof}
	Consider the geodesic triangle $\geodesicTriangle{p}{p'}{q}$ with some $q \in \cM$.
	Using \eqref{eq:Flat_isomorphism} and \eqref{eq3:theo:cosine.law} with $p_i = p'$, $p_{i \oplus 1} = p$ and $p_{i \ominus 1} = q$ we can say that
	\begin{equation*}
		\dual{[\logarithm{p} p']^{\flat}}{\logarithm{p} q} + \dual{[\logarithm{p'} p]^{\flat}}{\logarithm{p'} q}
		=
		\riemannian{\logarithm{p} p'}{\logarithm{p} q} + \riemannian{\logarithm{p'} p}{\logarithm{p'} q}
		\ge
		d^2(p,p')
		.
	\end{equation*}
	Hence, we obtain
	\begin{equation*}
		\dual{[\logarithm{p} p']^{\flat}}{\logarithm{p} q} - F(q)
		\ge
		d^2(p,p') + \dual{[-\logarithm{p'} p]^{\flat}}{\logarithm{p'} q} - F(q)
		,
	\end{equation*}
	Taking the supremum with respect to $q$ on both sides, we conclude the proof.
\end{proof}
We now establish a result regarding the properness of the conjugate function, thereby generalizing a result from \cite[Lem.~3.4]{BergmannHerzogSilvaLouzeiroTenbrinckVidalNunez:2021:1}.
\begin{proposition}
	Let $F \colon \cM \to \overline{\R}$.
	If $F^* \colon \cotangentBundle \to \overline{\R}$ is proper, then $F$ is also proper.
\end{proposition}
\begin{proof}
	Since $F^*$ is proper by assumption we have $\dom F^*\neq \emptyset$.
	Choose some $(p,\xi) \in \dom F^*$.
	Using \cref{def:Fenconj} we can say that
	\begin{equation*}
		+\infty > F^*(p,\xi)
		\ge
		\dual{\xi}{\logarithm{p} q} - F(q)
		\quad \text{for all } q \in \cM
		.
	\end{equation*}
	Since $-(-\infty) = +\infty$, we have that $F(q) \neq -\infty$ for all $q\in \cM$.
	Now, we will show that $\dom F\neq \emptyset$.
	Suppose, by contraposition, that $F(q) = +\infty$ for all $q \in \cM$.
	This would imply that $\dual{\xi}{\logarithm{p} q} - F(q) = -\infty$ for all $q \in \cM$ and, consequently, $F^*(p,\xi) = -\infty$, which contradicts the fact $(p,\xi) \in \dom F^*$.
	Therefore, $\dom F\neq \emptyset$ and proof is complete.
\end{proof}

Due to the relationship between \cref{def:Fenconj} and \cite[Def.~3.1]{BergmannHerzogSilvaLouzeiroTenbrinckVidalNunez:2021:1} mentioned in \cref{rem:previous-def}, the proof of the following result follows directly from Lem.~3.7. and Prop.~3.9 of \cite{BergmannHerzogSilvaLouzeiroTenbrinckVidalNunez:2021:1}.
Its proof will therefore be omitted.

\begin{proposition}\label{prop:properties}
	Let $F,G \colon \cM \to \overline{\R}$ be two proper functions and suppose that $\alpha \in \R$ and $\lambda > 0$.
	Then the following statements hold.
	\begin{enumerate}
		\item
			If $F(q) \le G(q)$ for all $q \in \cM$, then $F^*(p,\xi) \ge G^*(p,\xi)$ for all $(p,\xi) \in \cotangentBundle$.
		\item
			If $G(q) = F(q) + \alpha$ for all $q \in \cM$, then $G^*(p,\xi) = F^*(p,\xi) - \alpha$ for all $(p,\xi) \in \cotangentBundle$.
		\item
			If $G(q) = \lambda F(q)$ for all $q \in \cM$, then $G^*(p,\xi) = \lambda F^*(p,\frac{\xi}{\lambda})$ for all $(p,\xi) \in \cotangentBundle$.
		\item\label{item:properties_Fenchel-Young}
			The Fenchel--Young inequality holds, \ie, for all $(p,\xi) \in \cotangentBundle$ we have
			\begin{equation*}
				F(q) + F^*(p,\xi)
				\ge
				\dual{\xi}{\logarithm{p} q}
				\quad \text{for all } q \in \cM
				.
			\end{equation*}
	\end{enumerate}
\end{proposition}

Now we present a result that shows the \emph{partial} convexity of the Fenchel conjugate \wrt the second argument.
\begin{proposition}
	\label{proposition:Fconjugate_is_partially_convex}
	Let $F \colon \cM \to \eR$ be any function and $p \in \cM$.
	Then the function $F^*(p,\cdot) \colon \cotangent{p} \to \eR$ is convex.
\end{proposition}
\begin{proof}
	We can infer from \cref{def:Fenconj} that
	\begin{equation*}
		F^*(p,\xi)
		=
		\sup_{q \in \cM} \paren[big]\{\}{\dual{\xi}{\logarithm{p} q} - F(q)}
		\quad \text{for } \xi \in \cotangent{p}[\cM]
	\end{equation*}
	is the supremum over a family of affine functions in $\xi$ on the vector space $\cotangent{p}[\cM]$.
	Its convexity \wrt $\xi$ is therefore a standard result.
\end{proof}

\begin{remark}\label{rem:Subgrad.ts}
	Let $p \in \cM$ and suppose that $F^*(p,\cdot) \colon \cotangent{p} \to \overline{\R}$ is proper.
	The subdifferential of $F^*(p,\cdot)$ at $\xi \in \cotangent{p}$, denoted by $\partial_2 F^*(p,\xi)$, is the set
	\begin{equation*}
		\partial_2 F^*(p,\xi)
		=
		\setDef[auto]{X \in \tangent{p}}{F^*(p,\xi') \ge F^*(p,\xi) + \dual{\xi' - \xi}{X} \quad \text{for all } \xi' \in \cotangentBundle}
		.
	\end{equation*}
\end{remark}
In the following statement we give a characterization of this subdifferential in terms of the conjugate function.
This result is a generalization of \cref{thm:Classical_characterization_subdifferential} to the Riemannian context.
\begin{theorem}\label{theo:subgrad-conj}
	Let $F \colon \cM \to \overline{\R}$ be a proper convex function.
	Then $\xi \in \partial F(p)$ holds if and only if $F^*(p,\xi) = - F(p)$.
\end{theorem}
\begin{proof}
	First we consider the case $p \in \dom F$.
	Suppose that $\xi \in \partial F(p)$.
	Hence, using \cref{def:Subdifferential} we have
	\begin{equation*}
		\dual{\xi}{\logarithm{p} q} - F(q)
		\le
		- F(p)
		\quad \text{for all } q \in \cM
		.
	\end{equation*}
	Taking the supremum with respect to $q$ and considering \cref{def:Fenconj} we obtain
	\begin{equation*}
		F^*(p,\xi) = \sup_{q \in \cM} \paren[big]\{\}{\dual{\xi}{\logarithm{p} q} - F(q)}
		\le
		-F(p)
	\end{equation*}
	holds.
	On the other hand, using \cref{def:Fenconj} it is easy to see that
	\begin{equation*}
		F^*(p,\xi) = \sup_{q \in \cM} \paren[big]\{\}{\dual{\xi}{\logarithm{p} q} - F(q)}
		\ge
		\dual{\xi}{\logarithm{p} p} - F(p)
		=
		-F(p)
		.
	\end{equation*}
	Thus, $F^*(p,\xi) = - F(p)$ follows.

	For the converse, suppose that $\xi \in \cotangent{p}$ is chosen such that $F^*(p,\xi) = - F(p)$ holds.
	Hence, using \cref{def:Fenconj} we have
	\begin{equation*}
		- F(p)
		=
		F^*(p,\xi)
		\ge
		\dual{\xi}{\logarithm{p} q} - F(q)
		\quad \text{for all } q \in \cM
		.
	\end{equation*}
	Therefore, it follows from \cref{def:Subdifferential} that $\xi \in \partial F(p)$ holds.

	When $p \not \in \dom F$, then $F(p) = \infty$ and therefore $\partial F(p) = \emptyset$ since $F$ is proper.
	Suppose that there exists $\xi \in \cotangent{p}$ such that $F^*(p,\xi) = - F(p)$ holds.
	Proceeding as above this entails $- F(p) = F^*(p,\xi) \ge \dual{\xi}{\logarithm{p} q} - F(q)$ for all $q \in \cM$ and therefore $\xi \in \partial F(p)$, which is a contradiction.
	This concludes the proof.
\end{proof}

\begin{remark}
	In case $\cM = \R^n$, \cref{theo:subgrad-conj} reads:
	$\xi \in \partial F(p)$ if and only if
	\begin{equation*}
		-F(p)
		=
		F^*(p,\xi)
		=
		F^*(0,\xi) - \dual{\xi}{p}
		,
	\end{equation*}
	where the last equality follows from \eqref{eq:redundancy_in_the_conjugate_on_Rn}.
	Since $F^*(0,\xi)$ coincides with the classical definition of the Fenchel conjugate of $F$, we can indeed conclude that \cref{theo:subgrad-conj} generalizes \cref{thm:Classical_characterization_subdifferential} with $\cX = \R^n$ to the Riemannian case.
\end{remark}

The following result shows that, under certain conditions, a function $F \colon \cM \to \eR$ is bounded from below by a particular continuous function.
This function depends on the metric of $\cM$ and, in the Euclidean case, it is an affine function.
A version of this result in Euclidean spaces, whose proof we are following, can be found in \cite[Thm.~2.2.6]{Zalinescu:2002:1}.

\begin{lemma}\label{lem:aux.main.theo}
	Let $F \colon \cM \to \eR$ be a proper lsc convex function and $p \in \dom F$.
	Then there exist $q \in \dom F$, $\alpha \in \R$ and $\lambda > 0$ such that
	\begin{equation*}
		\lambda \, \riemannian{\logarithm{q} p}{\logarithm{q} q'} - F(q')
		\le
		\alpha
		\quad \text{for all } q' \in \dom F
		.
	\end{equation*}
\end{lemma}
\begin{proof}
	First take $s < F(p)$, \ie, $(p,s) \not \in \epi F$.
	Since $F$ is a proper lsc convex function, applying \cref{coro:ineq.proj} we conclude that there exist $(\hat p,\hat s) \in \epi F$ such that
	\begin{equation}\label{eq1:lem:aux.main.theo}
		\riemannian{\logarithm{\hat p} p}{\logarithm{\hat p} q'} + (s - \hat s)(r-\hat s)
		\le
		0
		\quad \text{for all } (q',r) \in \epi F
		.
	\end{equation}
	Taking $(q',r) = (p,F(p)+ n)$, $n \in \N$, we get
	\begin{equation}\label{eq2:lem:aux.main.theo}
		\riemannian{\logarithm{\hat p} p}{\logarithm{\hat p} p} + (s - \hat s)(F(p) + n - \hat s)
		\le
		0
		\quad \text{for all } n \in \N
		.
	\end{equation}
	From this $s - \hat s \neq 0$ follows, since otherwise we would have $s = \hat s$ and, by \eqref{eq2:lem:aux.main.theo}, $p = \hat p$ follows.
	Therefore we would have $(p,s) = (\hat p, \hat s)$, contradicting the fact $(p,s) \not \in \epi F$.
	On the other hand, considering \eqref{eq2:lem:aux.main.theo} with $s - \hat s > 0$ and $n$ sufficiently large, we obtain another contradiction.
	Therefore, we conclude $s - \hat s < 0$.
	Dividing \eqref{eq1:lem:aux.main.theo} by $\hat s - s > 0$, we have
	\begin{equation*}
		\frac{1}{\hat s - s} \riemannian{\logarithm{\hat p} p}{\logarithm{\hat p} q'} - r
		\le
		- \hat s
		\quad \text{for all } (q',r) \in \epi F
		.
	\end{equation*}
	Since $(q',F(q')) \in \epi F$ for all $q' \in \dom F$, it follows that
	\begin{equation*}
		\frac{1}{\hat s - s} \riemannian{\logarithm{\hat p} p}{\logarithm{\hat p} q'} - F(q')
		\le
		-\hat s
		\quad \text{for all } q' \in \dom F
		.
	\end{equation*}
	To finalize the proof choose $q = \hat p$, $\alpha = -\hat s$ and $\lambda = 1/(\hat s-s)$.
\end{proof}

The following result shows that our definition of Fenchel conjugate on $\cM$ allows us to obtain an extension of the second part of \cref{thm:Classical_Fenchel_Moreaou_theorem} to the Riemannian context.

\begin{theorem} \label{theo:F.conj.pro}
	Let $F \colon \cM \to \eR$ be a proper lsc convex function.
	Then $F^*$ is proper.
\end{theorem}
\begin{proof}
	Fix $(p,\xi) \in \cotangentBundle$ and choose some $p' \in \dom F$.
	Using \cref{rem:conjugate.dom}, we have
	\begin{equation*}
		F^*(p,\xi)
		=
		\sup_{q \in \dom F} \paren[big]\{\}{\dual{\xi}{\logarithm{p} q} - F(q)}
		\ge
		\dual{\xi}{\logarithm{p} p'} - F(p')
		>
		-\infty
		.
	\end{equation*}
	On the other hand, \cref{lem:aux.main.theo} guarantees that there are $q \in \dom F$, $\alpha \in \R$ and $\lambda > 0$ such that
	\begin{equation*}
		\lambda \, \riemannian{\logarithm{q} p'}{\logarithm{q} q'} - F(q')
		\le
		\alpha
		\quad \text{for all } q' \in \dom F
		.
	\end{equation*}
	Using \eqref{eq:Flat_isomorphism} and taking the supremum with respect to $q'$, we can conclude that $F^*\paren[big](){q,[\lambda \logarithm{q} p']^{\flat}} \le \alpha < +\infty$ holds.
	This shows $\dom F^*\neq \emptyset$, which completes the proof.
\end{proof}

\subsection{Fenchel Biconjugate}
\label{subsec:Fenchel-Biconjugate}

We now define the Fenchel biconjugate on Hadamard manifolds.

\begin{definition}\label{def:Fenconjbi}
	Let $F \colon \cM \to \eR$.
	The Fenchel biconjugate of $F$ is the function $F^{**} \colon \cM \to \eR$ defined by
	\begin{equation*}
		F^{**}(p)
		\coloneqq
		\sup_{(q,\xi) \in \cotangentBundle} \paren[big]\{\}{\dual{\xi}{\logarithm{q} p} - F^*(q,\xi) }
		\quad \text{for } p \in \cM
		.
	\end{equation*}
\end{definition}
Similarly as it was the case for $F^*$, our definition differs from our previous definition of the biconjugate in \cite[Def.~3.5]{BergmannHerzogSilvaLouzeiroTenbrinckVidalNunez:2021:1}.
In particular, $F^{**}$ does not depend on a base point.
The following remark shows that the above definition is a natural extension of \eqref{eq:Classical_biconjugate_function} to the Riemannian context.

\begin{remark}\label{rem:biconj.eucli.case}
	Let $F \colon \R^n \to \overline{\R}$.
	The \cref{def:Fenconjbi} with $\cM$ equal to the Euclidean space $\R^n$ becomes
	\begin{equation*}
		F^{**}(p)
		=
		\sup_{(q,\xi) \in \R^n \times \R^n} \paren[big]\{\}{\dual{\xi}{p-q} - F^*(q,\xi)}
		\quad \text{for } p \in \R^n
		.
	\end{equation*}
	Since \cref{rem:eucli.case}, \cref{eq:redundancy_in_the_conjugate_on_Rn} states that $F^*(q,\xi) = F^*(0,\xi) - \dual{\xi}{q}$ for all $(q,\xi) \in \R^n\times \R^n$, it follows that
	\begin{equation*}
		F^{**}(p)
		=
		\sup_{\xi \in \R^n} \paren[big]\{\}{\dual{\xi}{p} - F^*(0,\xi)}
		\quad \text{for } p \in \R^n
		.
	\end{equation*}
	Taking into account that $F^*(0,\xi)$ coincides with the classical definition of Fenchel conjugate on $\R^n$, we can conclude that \cref{def:Fenconjbi} generalizes the classical definition of the Fenchel biconjugate from the Euclidean space to Hadamard manifolds.
\end{remark}

The Fenchel biconjugate function is always a lower bound on the original function, as the following result states, which generalizes \cite[Prop.~13.14]{BauschkeCombettes:2011:1}.
\begin{proposition}\label{prop:biconj.less.fun}
	Let $F \colon \cM \to \overline{\R}$.
	Then $F^{**}(p) \le F(p)$ holds for all $p \in \cM$.
\end{proposition}
\begin{proof}
	Applying \cref{def:Fenconjbi} and \cref{def:Fenconj}, we have
	\begin{align*}
		F^{**}(p)
		&
		=
		\sup_{(q,\xi) \in \cotangentBundle} \paren[big]\{\}{\dual{\xi}{\logarithm{q} p} - F^*(q,\xi)}
		,
		\\
		&
		=
		\sup_{(q,\xi) \in \cotangentBundle} \paren[big]\{\}{\dual{\xi}{\logarithm{q} p} - \sup_{q' \in \cM} \paren[big]\{\}{\dual{\xi}{\logarithm{q} q'} - F(q')} }
		,
		\\
		&
		=
		\sup_{(q,\xi) \in \cotangentBundle} \paren[big]\{\}{\dual{\xi}{\logarithm{q} p} + \inf_{q' \in \cM} \paren[big]\{\}{- \dual{\xi}{\logarithm{q} q'} + F(q') } }
		,
		\\
		&
		\le
		\sup_{(q,\xi) \in \cotangentBundle} \paren[big]\{\}{\dual{\xi}{\logarithm{q} p} - \dual{\xi}{\logarithm{q} p} + F(p) }
		,
		\\
		&
		=
		F(p)
		\quad \text{for any } p \in \cM
		.
	\end{align*}
\end{proof}

The following result is a version of the famous Fenchel--Moreau theorem in the Riemannian case, compare \cref{thm:Classical_Fenchel_Moreaou_theorem}.
\begin{theorem} \label{theo:main}
	Let $F \colon \cM \to \eR$ be a proper lsc convex function.
	Then $F^{**} = F$ holds.
\end{theorem}
\begin{proof}
	The proof generally follows along the lines of the analog result in vector spaces; see, \eg, \cite[Thm.~13.32]{BauschkeCombettes:2011:1} and \cite[Thm.~2.3.3]{Zalinescu:2002:1}.
		However, it is worth pointing out that even in the case vector space setting, our definition of the Fenchel conjugate~$F^*(p,\xi)$ is more general due to the extra argument~$p$ (replacing~$0$), and therefore the proof requires adaption.
		This is all the more true for the Hadamard manifold setting, where some additional algebraic manipulations need to be rewritten, and terms need to be grouped differently.
		We therefore consider it justified to provide the generalized proof here.
		For the interested reader we point out where adaptations were necessary compared to the proof in \cite[Thm.~2.3.3]{Zalinescu:2002:1}.
	
	Let $p \in \cM$ be arbitrary.
	Choose some $s \in \R$ such that $s < F(p)$ holds, \ie, $(p,s) \not \in \epi F$.
	Since $F$ is a proper lsc convex function, we can apply \cref{coro:ineq.proj} to conclude that there exists $(\hat p,\hat s) \in \epi F$ such that
	\begin{equation}\label{eq1:theo:main}
		\riemannian{\logarithm{\hat p} p}{\logarithm{\hat p} q} + (s - \hat s)(r-\hat s)
		\le
		0
		\quad \text{for all } (q,r) \in \epi F
		.
	\end{equation}
	The inequality above is a rearrangement of eq.(2.32) in \cite{Zalinescu:2002:1}. 
	Considering \eqref{eq1:theo:main} with $(q,r) = (\hat p,F(\hat p)+ n)$, $n \in \N$, we obtain $(s - \hat s)(F(\hat p)+ n-\hat s) \le 0$ for all $n \in \N$.
	Since $\hat p \in \dom F$ holds, the assumption $s - \hat s > 0$ yields a contradiction for $n$ sufficiently large.
	Therefore, we must have $s - \hat s \le 0$.
	With this in mind, we will prove $s \le F^{**}(p)$.

	First, let us assume $s - \hat s < 0$.
	Dividing \eqref{eq1:theo:main} by $\hat s - s > 0$ we get
	\begin{equation*}
		\frac{1}{\hat s - s} \riemannian{\logarithm{\hat p} p}{\logarithm{\hat p} p'} - r
		\le
		-\hat s < - s
		\quad \text{for all } (p',r) \in \epi F
		.
	\end{equation*}
	Using \eqref{eq:Flat_isomorphism} and the expression above with $p' \in \dom F$ and $r = F(p')$, we have
	\begin{align*}
		\MoveEqLeft
		\dual{[(\logarithm{\hat p} p)/(\hat s - s)]^{\flat}}{\logarithm{\hat p} p'} - F(p')
		\\
		&
		=
		\frac{1}{\hat s - s}\riemannian{\logarithm{\hat p} p}{\logarithm{\hat p} p'} - F(p')
		<
		-s
		\quad \text{for all } p' \in \dom F
		.
	\end{align*}
    This inequality is a rearrangement of the inequality in line~13 in the proof by \cite{Zalinescu:2002:1}.
    Taking the supremum with respect to $p' \in \dom F$ and considering \cref{rem:conjugate.dom}, it follows that $F^*\paren[big](){\hat p,[(\logarithm{\hat p} p)/( \hat s - s )]^{\flat}} \le -s$.
	Taking into account that $\dual{[(\logarithm{\hat p} p)/(\hat s - s)]^{\flat}}{\logarithm{\hat p} p} \ge 0$, the last inequality and \cref{def:Fenconjbi} yield
	\begin{equation*}
		s
		\le
		\dual{[(\logarithm{\hat p} p)/(\hat s - s)]^{\flat}}{\logarithm{\hat p} p} - F^*\paren[big](){\hat p,[(\logarithm{\hat p} p)/( \hat s - s )]^{\flat}}
		\le
		F^{**}(p)
		.
	\end{equation*}
   The previous inequality corresponds to a manipulation of line~15 in \cite{Zalinescu:2002:1}.
	Now, let us prove that $s \le F^{**}(p)$ also holds when we assume $s - \hat s = 0$.
	In this case, \eqref{eq1:theo:main} becomes
	\begin{equation}\label{eq3:theo:main}
		\riemannian{\logarithm{\hat p} p}{\logarithm{\hat p} q}
		\le
		0
		\quad \text{for all } q \in \dom F
		.
	\end{equation}
	The previous inequality parallels the one in line~17 in \cite{Zalinescu:2002:1}. 
	Using \cref{lem:aux.main.theo} for $\hat p \in \dom F$, there exist $q' \in \dom F$, $\alpha \in \R$ and $\lambda > 0$ such that
	\begin{equation*}
		\lambda \, \riemannian{\logarithm{q'} \hat p}{\logarithm{q'} q} - F(q)
		\le
		\alpha
		\quad \text{for all } q \in \dom F
		.
	\end{equation*}
	On the other hand, it is easy to see that \eqref{eq3:theo:main} is equivalent to the inequality
	\begin{equation*}
		\lambda \, \riemannian{\logarithm{\hat p} q'}{\logarithm{\hat p} q}
		+
		\lambda \, \riemannian{\sigma \logarithm{\hat p} p - \logarithm{\hat p} q'}{\logarithm{\hat p} q}
		\le
		0
		\quad \text{for all } q \in \dom F, \; \sigma > 0
		.
	\end{equation*}
	Adding the last two inequalities we get
	\begin{multline}\label{eq5:theo:main}
		\lambda \, \riemannian{\logarithm{q'} \hat p}{\logarithm{q'} q}
		+
		\lambda \, \riemannian{\logarithm{\hat p} q'}{\logarithm{\hat p} q}
		\\
		+
		\lambda \, \riemannian{\sigma \logarithm{\hat p} p - \logarithm{\hat p} q'}{\logarithm{\hat p} q} - F(q)
		\le
		\alpha
		\quad \text{for all } q \in \dom F, \; \sigma > 0
		.
	\end{multline}
	Using \eqref{eq3:theo:cosine.law} for the geodesic triangle $\geodesicTriangle{\hat p}{q'}{q}$ for $q \in \cM$, we can conclude that $0 \le \riemannian{\logarithm{q'} \hat p}{\logarithm{q'} q} + \riemannian{\logarithm{\hat p} q'}{\logarithm{\hat p} q}$ holds for all $q \in \dom F$.
	Thus, \eqref{eq5:theo:main} yields
	\begin{equation*}
		\lambda \, \riemannian{\sigma \logarithm{\hat p} p - \logarithm{\hat p} q'}{\logarithm{\hat p} q} - F(q)
		\le
		\alpha
		\quad \text{for all } q \in \dom F, \; \sigma > 0
		.
	\end{equation*}
	The inequality above is an adaptation of line~19 in the proof by \cite{Zalinescu:2002:1}. 
	Considering \eqref{eq:Flat_isomorphism} and taking the supremum over $q \in \cM$, we get
	$F^*\paren[big](){\hat p, [\lambda \, (\sigma \logarithm{\hat p} p - \logarithm{\hat p} q')]^{\flat}} \le \alpha$ for all $\sigma > 0$.
	Therefore,
	\begin{align*}
		\MoveEqLeft
		\lambda \, \sigma \, d^2(p,\hat p) - \lambda \, \riemannian{\logarithm{\hat p} q'}{\logarithm{\hat p} p} - \alpha
		\\
		&
		=
		\riemannian{\lambda \, (\sigma \logarithm{\hat p} p - \logarithm{\hat p} q')}{\logarithm{\hat p} p} - \alpha
		\\
		&
		=
		\dual{[\lambda \, (\sigma \logarithm{\hat p} p - \logarithm{\hat p} q') ]^{\flat}}{\logarithm{\hat p} p} - \alpha
		\\
		&
		\le
		\dual{[\lambda \, (\sigma \logarithm{\hat p} p - \logarithm{\hat p} q')]^{\flat}}{\logarithm{\hat p} p} - F^*\paren[big](){\hat p, [\lambda \, (\sigma \logarithm{\hat p} p - \logarithm{\hat p} q')]^{\flat}}
		\\
		&
		\le
		\sup_{(q,\xi) \in \cotangentBundle} \paren[big]\{\}{\dual{\xi}{\logarithm{q} p} - F^*(q,\xi) }
		=
		F^{**}(p)
		,
	\end{align*}
	for all $\sigma > 0$. 
	The above chain of inequalites represent an adaptation of lines~21--23 in \cite{Zalinescu:2002:1}.
	As we are analyzing the case $s - \hat s = 0$, we can conclude that $p \neq \hat p$ must hold, otherwise we would have $(p,s) = (\hat p, \hat s)$, contradicting the fact $(p,s) \not \in \epi F$.
	Thus, taking $\sigma$ sufficiently large, we get $F^{**}(p) = +\infty$ and thus $s \le F^{**}(p)$ holds in this case as well.

	We have thus proved $s \le F^{**}(p)$ in all cases.
	Since $s < F(p)$ was arbitrary, we get $F(p) \le F^{**}(p)$.
	The conclusion of the proof now follows from \cref{prop:biconj.less.fun}.
\end{proof}

\subsection{Potential Application and Example}
\label{subsec:application}

In this section we briefly touch upon a potential application of the theory of Fenchel duality.
Although this is not explored further in the present paper, we mention that duality is a core concept that many solution algorithms for convex minimization problems are based on.
Specifically, suppose that $F \colon \cM \to \eR$ and $G \colon \cM \to \eR$ are proper, lsc and convex functions.
It follows from \cref{theo:main} and \cref{def:Fenconjbi} that the minimization problem 
\begin{equation*}\label{eq:primal.problem}
	\text{Minimize}
	\quad
	F(p) + G(p)
	,
	\quad
	p \in \cM
\end{equation*}
has the following saddle-point formulation:
\begin{equation*}
	\text{Minimize}
	\quad
	\sup_{(q,\xi) \in \cotangentBundle} \paren[big]\{\}{\dual{\xi}{\logarithm{q} p} + G(p) - F^*(q,\xi)}
	,
	\quad
	p \in \cM
	.
\end{equation*}
This formulation is the starting point for primal-dual algorithms, whose development on Hadamard manifolds is a topic for further research.
In any case, the evaluation of the conjugate $F^*$ of $F$ will be a requirement for the application of any such algorithm to a particular problem.

Therefore, we develop in this section a concrete example for the conjugate of a function.
$\cM$ will be the manifold $\spd{n}$ of real, symmetric, positive definite $n \times n$-matrices; see \cref{ex:nosingleton}. 
The function under consideration is
\begin{equation}
	\label{eq:log_det}
	F(B)
	\coloneqq
	a \, \ln \det(B)
\end{equation}
for some $a \in \R$.
This function appears in optimization problems from different fields.
We mention operator scaling, see \cite[eq.(1.1)]{AllenZhuGargLiOliveiraWigderson:2018:1}, as well as optimal experimental design, see \cite[eq.(5.8)]{PronzatoPazman:2013:1}, as two examples.
It also appears in minimum-volume covering problems, which have broad connections to many other fields; see \cite[Ch.~1.4]{Todd:2016:1}.

We will denote the identity matrix $n \times n$ by $\id$, and we will use the following property:
\begin{equation}\label{eq:pro:trln=lndet}
	\trace(\matrixLogarithm(C)) 
	= 
	\ln \det(C) 
	\quad
	\text{for all } 
	C \in \spd{n}
	.
\end{equation}
We recall that the map $\matrixLogarithm$ on the left-hand side of the above equality is the matrix logarithm, while $\ln$ on the right side is the natural logarithm of positive real numbers. 
The matrix exponential will be denoted by $\matrixExponential$. 

\begin{example}
	Consider the function $F \colon \cM \to \R$ from \eqref{eq:log_det} with some $a \in \R$.
	The manifold $\cM = \spd{n}$ is endowed with the Riemannian metric from \cref{ex:nosingleton}.
	Using \eqref{eq:spd_metric}--\eqref{eq:spd_logarithmic_map}, we evaluate
	\begin{align*}
		\dual{X^{\flat}}{\logarithm{A} B}  
		& 
		= 
		\trace \paren[big](){X A^{-1} [A^{1/2} \matrixLogarithm \paren[big](){A^{-1/2} B \, A^{-1/2}} \, A^{1/2}] A^{-1}} 
		\\
		& 
		= 
		\trace \paren[big](){X A^{-1/2} \matrixLogarithm \paren[big](){A^{-1/2} B \, A^{-1/2}} \, A^{-1/2}} 
		\\
		& 
		= 
		\trace \paren[big](){A^{-1/2}X A^{-1/2} \matrixLogarithm \paren[big](){A^{-1/2} B \, A^{-1/2}}}
	\end{align*}
	for any $A,B \in \cM$ and $X\in \tangent{A}$. 
	Using \cref{def:Fenconj} and the expression above, and performing the change of variable  $C = A^{-1/2} B \, A^{-1/2}$, we get
	\begin{align*}
		F^*(A,X^{\flat})
		&
		=
		\sup_{B \in \cM} \paren[big]\{\}{\dual{X^{\flat}}{\logarithm{A} B} - F(B)} 
		\\
		&
		=	
		\sup_{B \in \cM} \paren[big]\{\}{\trace \paren[big](){A^{-1/2}X A^{-1/2} \matrixLogarithm \paren[big](){A^{-1/2} B \, A^{-1/2}}} - a \ln \det(B)} 
		\\
		&
		=	
		\sup_{C \in \cM} \paren[big]\{\}{\trace \paren[big](){A^{-1/2}X A^{-1/2} \matrixLogarithm(C)} -  a \ln \det(A^{1/2} C \, A^{1/2})} 
		\\
		&
		=	
		\sup_{C \in \cM} \paren[big]\{\}{\trace \paren[big](){A^{-1/2}X A^{-1/2} \matrixLogarithm(C)} -  a \ln \det(C) -   a \ln \det(A)}
	\end{align*}
	for any $A \in \cM$ and $X\in \tangent{A}$. 
	Applying now \eqref{eq:pro:trln=lndet} and the linearity of $\trace$ and rearranging terms, we can obtain
	\begin{align*}
		F^*(A,X^{\flat})
		&
		=	
		\sup_{C \in \cM} \paren[big]\{\}{\trace \paren[big](){A^{-1/2}X A^{-1/2} \matrixLogarithm(C)} -  a\trace(\matrixLogarithm(C))    -   a\ln (\det( A  )) }
		\\
		&
		=	
		\sup_{C \in \cM} \paren[big]\{\}{\trace \paren[big](){A^{-1/2}X A^{-1/2} \matrixLogarithm(C)} -  \trace(a\matrixLogarithm(C))}  -   a\ln (\det( A  ))
		\\
		&
		=	
		\sup_{C \in \cM} \paren[big]\{\}{\trace \paren[big](){A^{-1/2}X A^{-1/2} \matrixLogarithm(C) -  a\matrixLogarithm(C)}}  -   a\ln (\det( A  ))
		\\
		&
		=	
		\sup_{C \in \cM} \paren[big]\{\}{\trace \paren[big](){\paren[big](){A^{-1/2}X A^{-1/2}   - a \, \id} \matrixLogarithm(C)}}  -   a\ln (\det( A  ))
		.
	\end{align*}
	The above calculation implies $F^*(A,X^{\flat}) = -a \ln \det(A)$ whenever $A^{-1/2}X A^{-1/2} = a \, \id$, \ie, $X = a \, A$.
	On the other hand, whenever $X \neq a \, A$, then choosing $C = \matrixExponential \paren[auto](){\lambda \, (A^{-1/2} X A^{-1/2} -  a \, \id)}$ for sufficiently large $\lambda > 0$ implies
	\begin{align*}
		\MoveEqLeft
		F^*(A,X^{\flat})
		\\
		&
		\geq 
		\sup_{\lambda > 0} \paren[big]\{\}{\trace \paren[big](){\paren[big](){A^{-1/2}X A^{-1/2} - a \, \id} \matrixLogarithm \paren[big][]{\matrixExponential \paren[big](){\lambda \, (A^{-1/2} X A^{-1/2} -  a \, \id)}}}} 
		\\
		&
		\quad
		- a  \ln \det(A) 
		\\ 
		&
		\geq 
		\sup_{\lambda > 0} \paren[big]\{\}{\lambda \trace \paren[big](){\paren[big](){A^{-1/2}X A^{-1/2} -  a \, \id}^2 }} - a \ln \det(A) 
		\\ 
		&
		=
		\sup_{\lambda > 0} \paren[big]\{\}{\lambda \, \norm{A^{-1/2}X A^{-1/2} - a \, \id}_F^2} - a \ln \det(A)
		\\
		&
		=
		\infty
		.
	\end{align*}
	Here $\norm{\cdot}_F$ denots the Frobenius norm.
	Overall, we conclude 
	\begin{equation*}
		F^*(A,X^{\flat})
		=
		\begin{cases}
			- a \ln \det(A) 
			& 
			\text{ if } 
			X = a \, A
			,
			\\
			+\infty 
			& 
			\text{ if } 
			X \neq a \, A
			.
		\end{cases}
	\end{equation*}
\end{example}
For comparison, we mention that the Fenchel conjugate of $F$ is different when we employ the classical conjugation concept from the ambient vector space $\symmetric{n}$ of symmetric $(n \times n)$-matrices, endowed with the Frobenius inner product.
In this case, $F^*(A) = \ln \det(-A)^{-1} - n = - \ln \det(-A) - n$ holds for $A \in -\spd{n}$ and $F^*(A) = \infty$ otherwise.
We refer the reader, \eg, to \cite[Ex.~3.23]{BoydVandenberghe:2004:1}.

\section{Separation of Convex Sets on Hadamard Manifolds}
\label{section:separation-theory}

Throughout this section, we develop a partial theory of separation of convex sets on a Hadamard manifold $\cM$ by affine hypersurfaces.
To see this theory on normed vector space, we refer the reader, \eg, to \cite[Ch.~1]{Brezis:2011:1}.
We begin by introducing a concept that generalizes the definition of an affine hyperplane to the Riemannian context.

\begin{definition}
	\label{definition:hypersurface}
	An \emph{affine hypersurface} of $\cM$ is a set $\cH \subset \cM$ of the form
	\begin{equation*}
		\cH
		=
		\hypersurface{p}{\xi}{\alpha}
		\coloneqq
		\setDef{q \in \cM }{\dual{\xi}{\logarithm{p} q} = \alpha}
		,
	\end{equation*}
	where $(p,\xi) \in \cotangentBundle$ and $\alpha \in \R$ are given with $\xi \neq 0$.
\end{definition}

\begin{remark}
	\label{rem:hypersurface-convex}
	Unlike in the Euclidean case, the affine hypersurface $\cH = \hypersurface{p}{\xi}{\alpha}$ is in general not totally geodesic, \ie, there might exist $q,q' \in \cH$ such that the unique geodesic segment $\geodesic<a>{q}{q'}$ in $\cM$ does not lie in $\cH$.
	This is due to the fact that $\cH$ is constructed via a hyperplane in the tangent space $\tangent{p}$, which only assures that all geodesic segments $\geodesic<a>{q}{p}$ are in $\cH$ if $\alpha = 0$.
	Otherwise, the only guarantee we have is that the connecting line $t\logarithm{p}q + (1-t)\logarithm{p}q'\in \tangent{p}$, $t \in \R$, satisfies $\dual{\xi}{t\logarithm{p}q + (1-t)\logarithm{p}q'} = \alpha$ and hence $c(t) = \exponential[big]{p}(t\logarithm{p}q + (1-t)\logarithm{p}q') \in \cH$.
	But this curve $c(t)$ is not necessarily a geodesic.
\end{remark}

\begin{remark}
	Consider the equivalence relation~$\sim$ defined in \eqref{eq:equivalence_relation_on_cotangent_bundle}.
	Note that $\hypersurface{p}{\xi}{\alpha} = \hypersurface{p'}{\xi'}{\alpha}$ holds for all $(p',\xi') \sim (p,\xi)$.
\end{remark}

\begin{definition}
	Let $\cA$ and $\cB$ be two subsets of $\cM$.
	We say that the hypersurface $\hypersurface{p}{\xi}{\alpha}$ \emph{separates} $\cA$ and $\cB$ if
	\begin{equation}
		\label{eq:separation_inequality}
		\dual{\xi}{\logarithm{p} q}
		\le
		\alpha
		\le
		\dual{\xi}{\logarithm{p} q'}
		\quad \text{holds for all } q \in \cA, \; q' \in \cB
		.
	\end{equation}
	We say that $\hypersurface{p}{\xi}{\alpha}$ \emph{strictly separates} $\cA$ and $\cB$ when both inequalities above are strict.
\end{definition}
Geometrically speaking, \eqref{eq:separation_inequality} means that $\cA$ lies in one of the \eqq{half-manifolds} determined by $\cH$, and $\cB$ lies in the other.

Let $\cA$ be a subset of $\cM$.
It is well known that the function $\phi \colon \cM \to \R$ defined by
\begin{equation}\label{eq:dist.fun.set}
	\phi(p)
	=
	d(p,\cA)
	\coloneqq
	\inf_{q \in \cA} d(p,q)
\end{equation}
is continuous on $\cM$.
This property will be used in the proof of the following result, which extends the classical strict separation theorem to the Riemannian setting.

\begin{theorem}\label{theo:separation.strict}
	Let $\cA \subset \cM$ and $\cB \subset \cM$ be two nonempty convex subsets such that $\cA \cap \cB = \emptyset$, $\cA$ is closed and $\cB$ is compact.
	Then there exists a hypersurface which strictly separates $\cA$ and $\cB$.
\end{theorem}
\begin{proof}
	Throughout the proof, points in $\cB$ will be marked by a prime.
	Since $\phi$ defined as in \eqref{eq:dist.fun.set} is continuous and $\cB$ is compact, the problem of minimizing $\phi$ over $\cB$ possesses at least one global solution.
	We denote one such solution by $\hat q' \in \cB$, \ie, $d(\hat q',\cA) \le d(q,\cA)$ holds for all $q \in \cB$.

	As $\cA$ is convex and closed, the projection map $\proj{\cA} \colon \cM \to \cA$ is well defined.
	Hence, setting $\hat q \coloneqq \proj{\cA}(\hat q')$, we have
	\begin{align*}
		d(q',\hat q)
		&
		\ge
		d(q', \cA )
		\ge
		\min_{q' \in \cB}d( q' , \cA )
		=
		\phi(\hat q')
		=
		d( \hat q' , \cA )
		\\
		&
		=
		d( \hat q' , \proj\cA(\hat q') )
		=
		d( \hat q', \hat q )
		\quad \text{for all } q' \in \cB
		,
	\end{align*}
	which means that $\hat q' = \proj{\cB}(\hat q)$.
	Taking into account $\cA \cap \cB = \emptyset$ we deduce $\hat q \neq \hat q'$.

	Let us define $p$ to be the midpoint of the geodesic segment connecting $\hat q$ to $\hat q'$.
	Then it is easy to see that we have
	\begin{equation}\label{eq2:theo:separation.strict}
		d(\hat q,p)
		=
		d(\hat q',p)
		=
		\frac{1}{2}
		d(\hat q,\hat q')
		>
		0
		.
	\end{equation}

	Next we prove $\proj{\cA}(p) = \hat q$ and $\proj{\cB}(p) = \hat q'$.
	Suppose by contradiction that $\proj{\cA}(p)\neq\hat q$ and consider the geodesic triangle $\geodesicTriangle{\hat q}{p}{\proj{\cA}(p)}$.
	Since $\logarithm{\hat q} p = \frac{1}{2} \logarithm{\hat q} \hat q'$ and $\hat q = \proj{\cA}(\hat q')$, \cref{theo:ineq.proj} guarantees
	\begin{align}
		\riemannian{\logarithm{\hat q} p}{\logarithm{\hat q} q}
		=
		\frac{1}{2}\riemannian{\logarithm{\hat q} \hat q'}{\logarithm{\hat q} q}
		&
		\le
		0
		\quad \text{for all } q \in \cA
		,
		\label{eq3:theo:separation.strict}
		\\
		\riemannian{\logarithm{\proj{\cA}(p)} p}{\logarithm{\proj{\cA}(p)} q}
		&
		\le
		0
		\quad \text{for all } q \in \cA
		.
		\label{eq4:theo:separation.strict}
	\end{align}
	Taking \eqref{eq3:theo:separation.strict} with $q = \proj{\cA}(p)$ and \eqref{eq4:theo:separation.strict} with $q = \hat q$ we get
	\begin{equation*}
		\riemannian{\logarithm{\hat q} p}{\logarithm{\hat q}\proj{\cA}(p)} + \riemannian{\logarithm{\proj{\cA}(p)} p}{\logarithm{\proj{\cA}(p)} \hat q}
		\le
		0
		,
	\end{equation*}
	which contradicts \eqref{eq3:theo:cosine.law} for $\geodesicTriangle{\hat q}{p}{\proj{\cA}(p)}$.
	Thus, we can conclude that $\proj{\cA}(p) = \hat q$ holds.
	Acting analogously with the geodesic triangle $\geodesicTriangle{\hat q'}{p}{\proj{\cB}(p)}$, we can also conclude that $\proj{\cB}(p) = \hat q'$.

	Consider the geodesic triangle $\geodesicTriangle{\hat q}{p}{q}$, $q \in \cA$.
	Since $\proj{\cA}(p) = \hat q$, \cref{theo:ineq.proj} and \eqref{eq3:theo:cosine.law} guarantee that
	\begin{equation*}
		\begin{alignedat}{2}
			-\riemannian{\logarithm{\hat q} p}{\logarithm{\hat q} q}
			&
			\ge
			0
			&
			&
			\quad \text{for all } q \in \cA
			,
			\\
			\riemannian{\logarithm{\hat q} p}{\logarithm{\hat q} q} + \riemannian{\logarithm{p} \hat q}{\logarithm{p} q}
			&
			\ge
			d^2(p,\hat q)
			&
			&
			\quad \text{for all } q \in \cA
			.
		\end{alignedat}
	\end{equation*}
	Adding the two inequalities above and using \eqref{eq2:theo:separation.strict} we can deduce that
	\begin{equation}\label{eq5:theo:separation.strict}
		\riemannian{\logarithm{p} \hat q}{\logarithm{p} q}
		\ge
		d^2(p,\hat q)
		>
		0
		\quad \text{for all } q \in \cA
		.
	\end{equation}
	Similarly, considering the geodesic triangle $\geodesicTriangle{\tilde{q}'}{p}{q'}$, $q' \in \cB$, and taking into account that $\proj{\cB}(p) = \hat q'$, we can also say that
	\begin{equation*}
		\riemannian{\logarithm{p} \hat q'}{\logarithm{p} q'}
		\ge
		d^2(p,\hat q')
		>
		0
		\quad \text{for all } q' \in \cB
		.
	\end{equation*}
	Since $\logarithm{p} \hat q' = -\logarithm{p} \hat q$, the last inequality implies $\riemannian{\logarithm{p} \hat q}{\logarithm{p} q'} < 0$ for all $q' \in \cB$.
	Hence, using \eqref{eq5:theo:separation.strict} and \eqref{eq:Flat_isomorphism} we can conclude that the hypersurface $\hypersurface[big]{p}{[\logarithm{p} \hat q]^{\flat}}{0}$ strictly separates $\cA$ and $\cB$.
\end{proof}

\section{Conclusions}
\label{sec:Conclusions}

In this paper we introduced a new definition of the Fenchel conjugate for functions defined on Hadamard manifolds.
In contrast to previous definitions, it is independent of the choice of a base point.
Our concept generalizes the Fenchel conjugate in the Euclidean case, and essential properties carry over.
As a next step we plan to investigate how to leverage the new concept algorithmically.
Moreover, we expect that a weaker version of the separation theorem can be shown, which merely requires $\cA$ and $\cB$ to be convex and closed.

\printbibliography

\end{document}